\documentstyle[12pt]{article}
\newtheorem{theorem}{THEOREM}[section]

\newtheorem{corollary}[theorem]{Corollary}
\newtheorem{lemma}[theorem]{Lemma}

\newtheorem{remark}{REMARK}

\newtheorem{proposition}[theorem]{Proposition}
\def\wbar{\overline{w}}

\def\lbar{\overline{\ell}}

\def\jbar{\overline{j}}

\def\xibar{\overline{\xi}}

\def\qbar{\overline{q}}
\def\hbar{\overline{h}}

\def\Omegabar{\overline{\Omega}}

\def\CC{{\rm\kern.24em\vrule
width.02em height1.4ex
depth-.05ex\kern-.26em C}}
\def\QQ{{\rm\kern.24em\vrule width.02em
height1.4ex depth-.05ex\kern-.26em Q}}
\def\RR{{\rm I\kern-.2em R}}
\def\HH{{\rm I\kern-.2em H}}
\def\ZZ{{\rm\kern.26em\vrule width.02em
height0.5ex depth0ex\kern.04em\vrule width.02em
height1.47ex depth-1ex\kern-.34em Z}}
\def\Ibb#1{{\rm I\kern-.23em#1}}
\def\Ib#1{{\rm I\kern-.25em#1}}
\def\k#1{\kern#1em}
\def\vb#1{\vrule width.02em height1.4ex depth-.05ex}
\def\NN{\Ibb N}
\def\11{{\rm\k{.45}\vb0\k{-.142}1}}

\def\Re{\mbox{\rm Re}\,}

\def\ibar{\overline{i}}
\def\epf{\hskip.2in\vrule width.4pt height6.65pt
depth.15pt\vrule
width2.5pt height6.65pt depth-6.25pt\hskip-2.5pt\vrule
width2.5pt
height.25pt depth.15pt\vrule width.4pt
height6.65pt depth.15pt\ }

\def\proof{\noindent {\bf Proof. } }

\def \Omegabar{\overline \Omega}

\def \d{\partial}

\def\zbar{\overline{z}}
\def\dbar{\overline{\partial}}
\def \hbar{\overline{h}}

\def\wbar{\overline{w}}

\def\\GBB{\cal B}

\def\1bar{\overline{1}}
\def\2bar{\overline{2}}

\def\dist{\hbox{dist}}
\def\Dbar{\overline{D}}
\def\qbar{\overline{q}}
\def\lambdabar{\overline{\lambda}}
\def\lambdabar{\overline{\lambda}}
\def\TP{\tilde{P}}
\begin{document}

 \title{ Solving  the Kerzman's problem on the 
 sup-norm  estimate for $\dbar$ on product domains}
\date{ Revised by July 10, 2023}
\author{Song-Ying Li}

\maketitle





					
	{\narrower\smallskip\noindent \noindent {\bf Abstract.}  In this paper, the author solves the long term open problem of Kerzman on sup-norm estimate for Cauchy-Riemann equation on polydisc
	in $n$-dimensional complex space. The problem has been open since 1971. He also extends and solves the problem on a bounded product domain $\Omega^n$,  where  $\Omega$ is any bounded domain in $\CC$  with $C^{1,\alpha}$ boundary for some $\alpha>0$.
	\par \smallskip}

\section{Introduction}

Let $\Omega$ be a bounded pseudoconvex domain in $\CC^n$. Let $f\in L^2_{(0,1)}(\Omega)$ be any $\dbar$-closed $(0,1)$-form
with coefficients $f_j\in L^2(\Omega)$. By H\"omander's theorem \cite{H1}, there is a unique $u\in L^2(\Omega)$ with $u\perp \hbox{Ker}(\dbar)$ such that $\dbar u=f$.
The regularity theory  for Cauchy-Riemann equations has been  a very important research area in several complex variables
for many decades. In particular, sup-norm estimate for $\dbar$
is the most difficult one. 
When $\Omega$ is a bounded smooth strictly pseudoconvex domain in $\CC^n$,
 in 1970,  Henkin \cite{Hen1},  Grauart and Lieb \cite{G-L} constructed a solution formula  for $\dbar u=f$  satisfying $\|u\|_{L^\infty }\le C_{\Omega}\|f\|_{L^\infty_{(0,1)}}$. 
 In 1971, Kerzman \cite{Ker} improved the above result in \cite{Hen1} and  \cite{GL}, 
 he  proved that $\|u\|_{C^\alpha (\Omega)}\le C_{\alpha, \Omega}\|f\|_{L^\infty_{(0,1)}}$
 for any $0<\alpha <1/2$. In 1971, Henkin and Romanov \cite{HR} proved the sharp estimate: $\|u\|_{C^{1/2}(\Omega)}\le C_\Omega \|f\|_{L^\infty_{(0,1)}}$.
 Recently, X. Gong \cite{Gon} generalized Henkin and Romanov's results.  He reduced  the assumption $\d\Omega\in C^\infty$ to $\d \Omega\in C^2$ and proved 
 that $\|u\|_{C^{\gamma +1/2}(\Omega)}
 \le \|f\|_{C^\gamma_{(0,1)}(\Omega)}$ for any $\gamma$ with that $\gamma+1/2$ is not an integer.
 In \cite{Ker}, when $\Omega=D^n$ is the unit polydisc in $\CC^n$,  Kerzman asked the following question: {\it  Does $\dbar u=f$
 have a solution satisfying $\|u\|_{C^\alpha }\le C_\alpha\|f\|_{L^\infty_{(0,1)}}$ for some $\alpha>0$ ?}  Let
 $f_j(\lambda)\in L^\infty(D) $ be  holomorphic in $D$ such that $u_0=\zbar_1 f_1(z_2)+\zbar_2 f_2(z_1) \not\in C(\overline{D}^2)$.  
 Let $f(z) =f_1(z_2) d\zbar_1+f_2(z_1) d\zbar_2$. Then $\dbar f=0$
 and $u_0\in L^\infty(D^2) \setminus C(\overline{D}^2)$ with $u_0\perp \hbox{Ker}(\dbar)$ solves $\dbar u=f$.
 Then the Kerzman's question can be refined by: {\it Does $\dbar u=f$ have a solution $u$ satisfying $\|u\|_{L^\infty}\le C\|f\|_{L^\infty_{(0,1)}} ?$} 
 The problem was studied by Henkin \cite{Hen}, he proved that if $f\in C^1_{(0,1)}(\Dbar^2)$ is $\dbar$-closed, then $\dbar u=f$ has a solution
 $u$ satisfying estimate $\|u\|_{L^\infty}\le C\|f\|_{L^\infty_{(0,1)}}$, where $C$ is a scalar constant. Notice that a $\dbar$-closed form
 $f \in L^\infty_{(0,1)}(D^n)$
 cannot  be approximated by $\dbar$-closed forms in $C^1_{(0,1)}(\overline{D}^n)$ in $L^\infty(D^n)$-norm. Henkin's result only partially answered Kerzman's 
 question and left  the Kerzman's question remanning open.

In \cite{Lan}, Landucci was able to improve the solution $u$ of $\dbar u=f$ in \cite{Hen} to the canonical solution which is the solution $u_0\perp \hbox{Ker}(\dbar)$. Recently,
Chen and McNeal \cite{CM} introduced a new  space ${\cal B}^p_{(0,1)}(D^n)$ of $(0,1)$-forms over $D^n$ which is  smaller than $L^p_{(0,1)}(D^n)$
 and proved $L^p$-norm estimates for $f\in {\cal B}^p_{(0,1)}(D^n)$ for $1<p\le \infty$. Their result
generalized Henkin's result. For a simple example, they reduced Henkin's assumption: $f=f_1 d\zbar_1+f_2 d\zbar_2 \in C^1_{(0,1)}(\overline{D}^2)$ to $f\in L^\infty_{(0,1)}(D^2)$
satisfying ${\d f_1 \over \d \zbar_2}\in L^\infty(D^2)$.  Dong, Pan and
Zhang \cite{DPZ} proved a very clean and pretty theorem: If $\Omega$ is any bounded domain in $\CC$ with $C^2 $ boundary
and  $f\in C_{(0,1)}(\Omegabar^n)$ is $\dbar$-closed, then the canonical solution $u_0$ of $\dbar u=f$ satisfies $\|u_0\|_{L^\infty}\le C\|f\|_{L^\infty_{(0,1)}}$.
However, $C_{(0,1)}(\Omegabar^n) $ is strictly smaller than $L^\infty_{(0,1)}(\Omega^n)$, the Kerzman's question remains open (see \cite{Li}).   

Main purpose of the current paper is to give a complete solution of the Kerzman's long open problem on 
the unit polydisc in $\CC^n$. More general, we will prove that the canonical solution $u$ satisfying estimate $\|u\|_{\infty}\le C\|f\|_\infty$
 on  the product domains $\Omega^n$, where $\Omega$ is a bounded domains $\Omega\subset \CC$ with $C^{1,\alpha}$ boundary for some $\alpha>0$.
   The main theorem is stated as follows.

\begin{theorem} Let $\Omega$ be a bounded domain in $\CC$ with $C^{1,\alpha}$ boundary for some $\alpha>0$.  Let
$f\in L^\infty_{(0,1)}(\Omega^n)$ be $\dbar$-closed. Then the canonical solution $u_0$ of  $\dbar u=f$ is constructed and satisfies
$$
\|u_0\|_{L^\infty(\Omega^n)}\le C\|f\|_{L^\infty_{(0,1)}(\Omega^n)}.\leqno(1.1)
$$
\end{theorem}

More information for $\dbar$-estimates, one may find  from the following references as well as the references therein. For examples,  Chen and Shaw \cite{CS}, Fornaess and Sibony \cite{FS}, 
Krantz \cite{Kran, K76}, Range \cite{Ran}, Range and Siu \cite{RS72, RS73}, Shaw \cite{S91}
and Siu \cite{Siu}. For product domains, one may also see  \cite{CS},   \cite{DLT},  \cite{KL2} and other related articles in the reference.

The paper is organized as follows. In section 2,  we first provide an estimate for the Green's function  and its derivatives on a bounded domain 
in $\CC$ with $C^{1,\alpha}$ boundary for some $\alpha>0$.  In Section 3, we provide  a formula solution for canonical solution of $\dbar u=f$  on the product domains.
In Section 4, technically, we translate the formula in Section 3 to one, from which we can get a uniform $L^p$ estimates. In Section 5,  we will prove Theorem 1.1. In fact,
we have proved  that Theorem 1.1 remains true for all $1<p\le \infty$.
Finally, in Section 6, based on $\dbar$ -estimate  on the disc $D\subset \CC$, we give a  more sharp theorem (Theorem 6.1) which is better  than Theorem 1.1.
\medskip

{\bf Acknowledgment.} The author would like to thank R-Y. Chen who read  through the first draft of manuscript and  Sun-sig Byun for providing some useful reference on Green's function (Theorem 2.1).
The author deeply appreciate Xianghong Gong  who went through the draft of the revision very carefully  and gave many invaluable comments and suggestions.
 The author also appreciates the referees for their invaluable comments and suggestions which are very helpful for the revision. 

\section{Green's function and Bergman kernel}

\subsection{Green's functions}

Let $\Omega$ be a bounded domain in $\CC$ and
let $G(\lambda,\xi)$ be the Green's function for the Laplace operator $\Delta'={\d^2 \over \d z \d \zbar}={1 \over 4} \Delta $ on 
$\Omega$. Then the Green's operator $G$ is defined by
$$
G[f](z)=\int_\Omega G(z, w) f(w) dA(w)\leqno(2.1)
$$
and  $G[f]$ satisfies
$$
{\d^2 G[f] \over \d \lambda \d\lambdabar}(\lambda) =f(\lambda),\lambda \in \Omega \quad \hbox{and}\quad G[f]=0 \ \hbox{ on } \d \Omega. \leqno(2.2)
$$
Notice that
$$\label{G1}
G(z,w)={1\over  \pi}  \log |z-w|^2-U(z, w), \leqno(2.3)
$$
where $U(z, w)$ is harmonic in $z$ and in $w$ separately which is the solution of
$$
\cases{ \Delta'_w U(z,w)=0,\quad w\in \Omega,\cr U(z,w)={1\over \pi} \log |z-w|^2,\quad w\in \d \Omega.\cr}\leqno(2.4)
$$ 
Then
$$
(2.5)\quad {\d ^2 G(z, w)\over \d z\d \wbar}=\delta_{w=z} -{\d^2 U(z, y)\over \d z \d \wbar},\  \hbox{and } \ {\d U(z,w)\over \d z}={1\over \pi(z-w)},\quad w\in \d\Omega.
$$
where $\delta_{w=z}=\delta_z(w)$ is the Dirac mass concentrated at $z$ as measure of $w$. Moreover, 
  $-{\d^2 U(z,w)\over \d z\d \wbar}$ is holomorphic in $z$ and anti-holomorphic in $w$. In fact, one can show that it is the Bergman kernel function for $\Omega$
  (see the details in \cite{Gar}).

Let $ A^2(\Omega)$ be the Bergman space over $\Omega$ which is the holomorphic subspace of $L^2(\Omega)$.
Let ${\cal P}:L^2(\Omega)\to A^2(\Omega)$ be the Bergman projection. Then
$$
(I-{\cal P})f (z)=-\int_\Omega {\d^2 G(z, w) \over \d z \d \wbar} f(w) dA(w),
$$
where the Bergman kernel function of $\Omega$ is given by
$$\label{BK1}
K(z,w)={\d^2 G(z,w)\over \d z\d \wbar}=-{\d^2 U(z,w)\over \d z\d \wbar}, \quad z\ne w.\leqno(2.6)
$$
By Theorem 0.5 in Jerison and Kenig \cite{JK}, {\it if $\d \Omega$ is Lipschitz, there is a $p_1>4$ such that the Green's operator $G: W^{-1, p}(\Omega)\to W^{1, p}(\Omega)$
is bounded for $p_1'<p<p_1$.} (2.6) and the expression of $I-{\cal P}$ imply that if $\d \Omega$ is Lipschitz, then ${\cal P}: L^p(\Omega)\to A^p(\Omega)$ is bounded for $p_1'<p<p_1$.
One may find further information on regularity of Bergman projections  in \cite{LL}.

 We need some properties of the Green's function and estimations on the Green's function and its derivatives  based on the regularity
 of $\d \Omega$. We recall  a definition. {\it We say that  a bounded domain $\Omega \subset \RR^n$  satisfies  a uniform exterior ball condition if there is a positive number $r$ such that
 for any $z_0\in \d \Omega$, there is $z_0\in \RR^n\setminus \Omegabar$ such that $\overline{B(z_0(r), r)} \cap \Omegabar=\{z_0\}$, where $B(x, r)$ is ball in $\RR^n$
  centered at $x$ with radius $r$.}  It is easy to see that if $\d \Omega$ is $C^2$, then
 $\Omega$ satisfies a uniform  exterior (and interior) ball condition.
 
The following theorem on the Green's function  was proved by  Gr\"uter and Widman  \cite{GW} (Theorem 3.3 for $n>2$).

\begin{theorem} If $\Omega$ is a bounded domain  in $\RR^n$ which satisfies a uniform exterior ball condition, then its associated Green function  satisfies the
following five properties for all $x, y\in \Omega$:

(i) $|G(x, y)|\le  C d_\Omega( x) |x-y|^{1-n}$;

(ii) $|G(x,y)|\le C d_\Omega (x)  d_\Omega (y)  |x-y|^{-n}$;

(iii) $ |\nabla_x G(x,y)|\le C|x-y|^{1-n}$;

(iv) $|\nabla_x G(x, y)| \le C d_\Omega (y) |x-y|^{-n}$;

(v) $|\nabla_x \nabla_y G(x, y)|\le C|x-y|^{-n} $.

\noindent Here $C$ is a constant depending only on $\Omega$ and $d_\Omega(x)$ is distance from $x$ to $\d \Omega$.
\end{theorem}


Notice that  $\Omega $ having $ C^{1,\alpha}$ boundary with $\alpha \in (0,1)$ may not  satisfy  a uniform exterior ball condition.
 We will give a formula for the Green's function on a bounded domain  in $\CC$ with $C^{1,\alpha}$ boundary. First, we give the formula
 for simply connected domain.

Applying the argument by Kerzman \cite{Ker2} and regularity theorem (Theorem 8.34 in \cite{GT}), one can prove the following result. 

\begin{proposition}\label{R1} Let $\Omega$ be a bounded domains in $\CC$ with
$C^{1,\alpha}$ boundary for some $0<\alpha<1$. Let $D(0,1)$ be the unit disc centered at $0$.

(i)  If $\psi: \Omega\to D(0,1)$ is a proper
holomorphic map, then $\psi\in C^{1,\alpha}(\Omegabar)$;

(ii) If $\phi: \Omega\to D(0,1)$ is biholomorphic, then the Green's function $G_\Omega$  for ${\d^2 \over \d z\d \zbar}$ in $\Omega$ is given by
$$\label{G2}
G_{\Omega}(z,w)={1\over  \pi} \log \Big|{\phi(z)-\phi(w)\over 1-\phi(z)\overline{\phi(w)}}\Big|^2 \leqno(2.7)
$$
which satisfies (i)--(v) in Theorem 2.1.
\end{proposition}

\proof By Theorem 8.34 in \cite{GT}, if $g\in L^\infty(D)$ with $D=D(0,1)$,  then 
$$
\Delta' u=g \hbox{ in } D, \quad u=0 \hbox{ on }\d D
$$
has a unique solution $u\in C^{1,\alpha}(\overline{D} )$. Let $g\in C^\infty_0 (D)$ be  a non-negative function on $D$ such that
 $\{z\in D: g(z)>0\}$
is a non-empty, relatively compact subset in $D$. Let $v(z)=u(\psi(z))$ be a function on $\Omega$ which solves the Dirichlet boundary problem:
$$
\cases{\Delta' v(z)=g(\psi(z))|\psi'(z)|^2,  \ z\in \Omega, \cr
\quad v(z)=0, \qquad\qquad \quad\   z\in \d \Omega.\cr}
$$
By the elliptic theory (Theorem 3.34 in \cite{GT}), one has $v\in C^{1,\alpha}(\Omegabar)$. Then
$$
{\d v\over \d z}(z)={\d u\over \d w}(\psi(z))\psi'(z).
$$
Since $D$ satisfies an  interior ball condition, by Hopf's lemma, one has ${\d u\over \d w}(w)\ne 0$ on $\d D$. Since
$u\in C^{1,\alpha}(\overline{D})$, one has ${\d u\over \d w}(w)\ne 0$ on the closed annulus
$A(0, 1-\epsilon, 1]=
\{w\in D: 1-\epsilon \le |w| \le 1\}$
for some small $\epsilon>0$. This implies
$$\label{Ph}
\psi'(z) ={\d v(z) \over \d z} /{\d u\over \d w}(\psi(z))\quad \hbox{ on } \psi^{-1}(A(0, 1-\epsilon, 1]). \leqno(2.8)
$$
 This implies that $\psi \in C^1(\Omegabar)$ since $\psi$ is holomorphic in $\Omega$. Applying (2.8) again, one can see that 
$\psi'(z)\in C^\alpha (\Omegabar)$. Therefore, $\psi \in C^{1,\alpha}(\Omegabar)$.

\medskip
It is well known that the Green's function for ${\d^2\over \d z \d \zbar}$ in the unit disc $D$ is:
$$
G(z,w)={1\over  \pi} \log\Big |{w-z\over 1-z\wbar}\Big |^2,\quad z, w\in D.\leqno(2.9)
$$
If $\phi: \Omega\to D$ is a bilomorphic map, then it is easy to check that the Greens's function for $\Omega$ is given by (2.7).
Moreover, one can check
that $G_\Omega$ satisfies Properties (i)--(v) in Theorem 2.1 when $n=2$. \epf
\medskip

When $\Omega$ is finite multiple connected, the Green's function of $\Omega$ was studied by S. Bell in \cite{Bel} by using Ahlfors map. We will use the partition of
unity  and (2.7) to study the Green's function. 

Let $\TP_\Omega [f]$ denote the  unique solution of  the Dirichlet boundary problem:
$$
\cases{\Delta'_w u(w)=0, \quad w\in \Omega ,\cr u( w)=f(w), \quad w\in \d \Omega.\cr}\leqno(2.10)
$$
Let $D=D(0,1)$ and let
$$
 r(w)=1-|w|^2, \quad r(z,w)=\cases{1-\langle z, w\rangle, z\in \Dbar,\cr z-w, \quad \hbox{if } z\not\in \Dbar},\quad w\in D.\leqno(2.11)
$$
Then for $f\in C^1(\Dbar)$
\begin{eqnarray*}
(2.12)\qquad \TP_D[f](w)&=&{1\over 2\pi}\int_{\d D}{r(w)\over | r(w,\xi)|^2} f(\xi) d\sigma(\xi)\\
&=&{r(w) \over \pi}\int_{D} {\d\over \d \xibar}\Big( {f(\xi) \xibar \over| r(w,\xi)|^2} \Big) A(\xi)\\
&=&{r(w) \over \pi}\int_{D} {\d\over \d \xi}\Big( {f(\xi) \xi \over |r(w,\xi)|^2} \Big) A(\xi),\quad w\in D.\qquad \qquad\qquad
\end{eqnarray*}
\begin{lemma} \label{Lem} For
$\psi \in C^{1,\alpha}(\Dbar)$ for some $\alpha>0$ with $D=D(0,1)$, we let
$$
H(z,w)=\TP_D [ \psi (\cdot) {\log|z-\cdot |^2\over \pi}](w)-\psi(w) \TP_D [{\log|z-\cdot|^2\over \pi}](w). \leqno(2.13)
$$
Then  for any $0<\epsilon\le 1$, there is constant $C_\epsilon$ depending only on $\epsilon$, $R$ and $\|\psi\|_{C^{1,\alpha}(\Dbar)}$ such that
for any $(z,w)\in D(0,R)\times D$ for any $R\ge 1$, one has
$$
|H(z,w)|\le C_\epsilon { r(w)  \over |r(z,w)|^\epsilon } \log \Big( {1+ r(w) \over r(w)}\Big),  \leqno(2.14)
$$
$$
|\nabla_z H(z, w)|\le {C_\epsilon r(w) \over |r(z,w)|^{1+\epsilon }} \leqno(2.15)
$$
 and
$$
\Big|\nabla_w {\d \over \d z} H(z,w)\Big|\le 
{C_\epsilon \over |r(z,w)|^{1+\epsilon} } .\leqno(2.16)
$$
\end{lemma}

\proof Define
$$
 X(w,\xi)= {\d \over \d \xibar} \Big( \xibar {\psi(\xi)-\psi(w) \over |r(w,\xi)|^2}\Big) 
\ \hbox{ and }\ 
Y(w,\xi)= {\psi(\xi)-\psi(w) \over |r(w,\xi)|^2}\xibar.\leqno(2.17)
$$
By (2.12),  (2.13) and (2.17), one has
\begin{eqnarray*}
(2.18)\quad H(z, w)
&=&{-i r(w) \over 2\pi^2} \int_{D} \Big[X(w,\xi)  \log |z-\xi|^2 +Y(w,\xi) {1\over \xibar-\zbar} \Big]
d\xibar\wedge d\xi \\
&=&{r(w) \over \pi^2}\Big( \int_D X(w,\xi) \log |z-\xi|^2 
d A(\xi) + \int_{D} { Y(w,\xi) \over \xibar-\zbar} 
dA(\xi )\Big).
\end{eqnarray*}
Coupling this with
$$
{\d \over \d \zbar}\Big( {1\over \pi}\int_D{g(w) \over z-w} dA(w)\Big)=g(z),
$$ 
one has
$$
{\d H(z,w)\over \d z}={r(w)\over \pi^2}\Big( \int_D X(w,\xi) {1\over z-\xi}
d A(\xi) + \pi \overline{Y(w,z)} \Big).\leqno(2.19)
$$
Similarly, if one let
$$
\tilde{X} (w,\xi)={\d \over \d \xi} \Big[\xi {\psi(\xi)-\psi(w) \over |r(w,\xi)|^2}\Big] \
\ 
\hbox{ and }\ 
\tilde{Y}(w,\xi)= {\psi_j(\xi)-\psi_j(w) \over |r(w,\xi)|^2}\xi.
$$
 Then
$$
{\d H(z,w)\over \d \zbar}={r(w) \over \pi^2}\Big( \int_D  \tilde{X}(w,\xi) {1\over\zbar-\xibar}
d A(\xi) + \pi \tilde{Y}(w,z) \Big).\leqno(2.20)
$$
Since  $ \psi\in C^1(\Dbar)$, one has
$$
 |X(w,\xi)|\le { C \over |r(w,\xi)|^2}\quad\hbox{and} \quad
|Y(w,\xi)|\le {C \over| r (w,\xi)|}.\leqno(2.21)
$$
where $C=C_0\|\psi\|_{C^1(\Dbar)}$, and
$$
 |\nabla_w X(w,\xi)|\le { C \over |r(w,\xi)|^3}\quad\hbox{and} \quad
|\nabla_w Y(w,\xi)|\le {C \over |r (w,\xi)|^2}.\leqno(2.22)
$$
For $0<\epsilon\le 1$, notice that  
$$
1={r(z,w)\over r(z,w)} ={r(\xi, w)+ \wbar (\xi-z) \over r(z, w)},\quad 1 \le {|r(w,\xi)|^\epsilon +|\xi-z|^\epsilon \over |r(z,w)|^\epsilon}\leqno(2.23)
$$
and
$$
\int_D {1  \over |1-w \xibar |^{\beta}} dA(\xi)\le
\cases{C \log ({1+ r(w) \over r(w)}),\quad \beta=2;\cr  {C\over 2 -\beta},\qquad \qquad \quad \beta <2; \cr 
{C\over \beta-2} r (w)^{2-\beta} ,\quad \ \  \beta>2.\cr}\leqno(2.24)
$$
For $\epsilon>0$, let
$$
C_\epsilon= \sup\{| t^\epsilon \log  t|: 0<t\le R\} +\| \log |z-\cdot|^2\|_{L^{2/\epsilon}(D)}<\infty.
$$
Then  $C_\epsilon\ge {1\over e \epsilon}$.  By (2.23),  (2.24) and the H\"older's inequality, one has 
\begin{eqnarray*}
\lefteqn{\int_D |X(w,\xi)| |\log |z-\xi|| dA(\xi)}\\
&\le&{1\over | r(z,w)|^\epsilon }\int_D  |X(w,\xi)| \Big( |r(w,\xi)|^\epsilon+|z-\xi|^\epsilon\Big)||\log|z-\xi|| dA(\xi)\\
&\le&{C\over | r(z,w)|^\epsilon }\int_D |r(w, \xi)|^{-2+\epsilon}  \, |\log|z-\xi|| dA(\xi)\\
&&+{C \over |r(z,w)|^\epsilon  }\int_D  |r (w,\xi)|^{-2} |z-\xi|^\epsilon  |\log|z-\xi||dA(\xi)\\
&\le &{C\over |r(z,w)|^\epsilon} \|{1\over |r(w,\cdot)|^{2-\epsilon}}\|_{L^{p'}}\|\log |z-\cdot |^2 \|_{L^p}+ {C_\epsilon\over |r(z,w)|^\epsilon} \log( {1+r(w)\over r(w)})\\
&\le& {C_\epsilon\over |r(z,w)|^\epsilon} \log( {1+r(w)\over r(w)})
\end{eqnarray*}
by taking $p'={2\over 2-\epsilon}$ and $p=2/\epsilon$. Since
\begin{eqnarray*}
\lefteqn{(2.25)\  \int_D { |r(w,\xi) |^\epsilon +|z-\xi|^\epsilon \over |r(w,\xi)| |z-\xi|}dA(\xi) }\\
&\quad\le& {C\over |r(z,w)|^\epsilon} \Big(\|{1\over r(w,\cdot)|^{1-\epsilon}}\|_{L^{2-\epsilon\over 1-\epsilon}}\|{1\over z-\xi}\|_{L^{2-\epsilon}}
+\|{1\over |z-\cdot |^{1-\epsilon}}\|_{L^{2-\epsilon\over 1-\epsilon}}\|{1\over  r(w, \cdot)}\|_{L^{2-\epsilon}}\Big)\\
&\quad\le&{C\over |r(z,w)|^\epsilon } \Big(\int_D |r(w,\xi)|^{-2+\epsilon} dA(\xi) \Big)^{1-\epsilon \over 2-\epsilon}\Big(\int_D |z-\xi|^{-2+\epsilon } dA(\xi)\Big)^{1 \over 2-\epsilon} \\
&&+{C\over |r(z,w)|^\epsilon } \Big(\int_D |r(w,\xi)|^{-2+\epsilon} dA(\xi) \Big)^{1 \over 2-\epsilon}\Big(\int_D |z-\xi|^{-2+\epsilon } dA(\xi)\Big)^{1-\epsilon \over 2-\epsilon} \\
&\quad \le &{C_\epsilon \over |r(z,w)|^\epsilon}.
\end{eqnarray*}
By (2.23) and (2.25), one has
\begin{eqnarray*}
\int_D |Y(w,\xi)|{1\over |z-\xi|} dA(\xi)
&\le&\int_D {C \over |r(w,\xi)|}{1\over |z-\xi|}dA(\xi) \\
&\le&{C\over |r(z,w)|^\epsilon } \int_D { |r(w,\xi) |^\epsilon +|z-\xi|^\epsilon \over |r(w,\xi)| |z-\xi|}dA(\xi) \\
&\le &{C_\epsilon \over |r(z,w)|^\epsilon}.
\end{eqnarray*}
Combining the above two estimates and (2.18) (the expression of $H$), one has proved (2.14). 
\medskip

By (2.23), one has
$$
\int_D  {1\over z-\xi} X(w, \xi)dA(\xi)
={1\over r(z,w)} \int_D{r(w,\xi) \over z-\xi} X(w,\xi) -\wbar X(w, \xi) dA(\xi).
$$
Since $(\xi-w)$ and $\xibar-\wbar$ is harmonic, $\xibar d\xi=i d\theta$ if $\xi=e^{i \theta}$, by Poisson integral formula, one has
$$
\int_{\d D}  {\xibar (\xi-w)\over |r(w, \xi)|^2} d\xi={2\pi \over i r(w) } (w-w)=0, \quad \int_{\d D}  {\xibar (\xibar-\wbar)\over |r(w, \xi)|^2} d\xi={2\pi \over i r(w)}(\wbar-\wbar)=0.
$$
Therefore,
\begin{eqnarray*}
\lefteqn{\Big|\int_D X(w, \xi)dA(\xi)\Big|}\\
&=&\Big|{1\over 2i} \int_{\d D} \xibar {\psi(\xi)-\psi(w)\over |r(w, \xi)|^2} d\xi\Big|\\
&\le& {1\over 2} \int_{\d D} |  {\psi(\xi)-\psi(w)-{\d\psi (w)\over \d w}(\xi-w)-{\d\psi(w)\over \d \wbar}(\xibar-\wbar)|\over |r(w, \xi)|^2} |\xibar d\xi|\\
&&+ {1\over 2} |{\d\psi (w)\over \d w}| |\int_{\d D}  {\xibar (\xi-w)\over |r(w, \xi)|^2} d\xi|
+ {1\over 2} |{\d\psi (w)\over \d \wbar}\int_{\d D}  {\xibar (\xibar-\wbar)\over |r(w, \xi)|^2} d\xi|\\
&\le& {C\|\psi\|_{C^{1,\alpha}}\over \alpha}.
\end{eqnarray*}

Applying (2.23)  with $0<\epsilon<1  $ and by (2.25),  one has
\begin{eqnarray*}
\lefteqn{\int_ D |X(w, \xi)|{|r(w,\xi)|\over |z-\xi|} dA(\xi)}\\
&\le&{C\over |r(z,w)|^{\epsilon}} \Big(\int_D  ( {|r(w, \xi)|^{\epsilon} \over |z-\xi|} +{1 \over |z-\xi|^{1-\epsilon} } )|r(w,\xi)|X(w,\xi)| dA(\xi)\\
&\le&{C\over |r(z,w)|^{\epsilon}} \int_D   {|r(w,\xi)|^\epsilon +|z-\xi|^\epsilon \over  |r(w, \xi)|  |z-\xi|} dA(\xi)\\
&\le &{ C_\epsilon \over |r(z,w)|^{\epsilon}} .
\end{eqnarray*}
By (2.21), one has that
$$
|Y(w, z)|\le {C\over |r(z,w)|}.
$$
Combing the above two estimates and (2.18), one has
$$
|{\d  H(z,w)\over \d z} |\le { C_\epsilon r(w) \over |r(z,w)|^{1+\epsilon}}.
$$
Using $\tilde{X}$ and $\tilde{Y}$ and  (2.20), one has the same estimates for $|{\d H(z, w)\over \d \zbar}|$. Therefore, (2.15)  is proved. 

By (2.25), one has
$$
\int_D {1\over |r(w,\xi)|^{2-\epsilon} |\xi-z|} dA(\xi)
\le {C\over r(w)} \int_D { |r(w,\xi)|^{\epsilon}\over |r(w,\xi)||z-\xi|} dA(w)
\le {C\over \epsilon r(w)}.
$$
Thus, applying the inequality (2.23) with $\epsilon=1$ first then with $0<\epsilon<1$, and the estimate for $\nabla_wX(w,\xi)$ in (2.22),  (2.24) and (2.25), one has
\begin{eqnarray*}
\lefteqn{\int_D \Big| {1\over |z-\xi|} |\nabla_wX(w, \xi)|\Big| dA(\xi)}\\
&\le&{C\over |r(z,w)|} \int_D   {|r(w, \xi)|\over |z-\xi|}|\nabla_wX(w,\xi)| dA(\xi)
+{C\over |r(z,w)|} \int_D |\nabla_wX(w,\xi)| dA(\xi)\\
&\le&{C\over |r(z, w)|^{1+\epsilon}  } \int_D ({r(w, \xi)|^{1+\epsilon} \over |z-\xi|} +{|r(w, \xi)| |z-\xi|^\epsilon\over |z-\xi|} ){1\over |r(w,\xi)|^3}  dA(\xi)+{C\over |r(z,w)| r(w)}\\
&\le&{C\over |r(z, w)|^{1+\epsilon} r(w) } \int_D {|r(w,\xi)|^\epsilon+|z-\xi|^\epsilon  \over |z-\xi| |r(w,\xi)|} dA(\xi)+{C  \over |r(z,w)|  r(w)}\\
&\le&{C_\epsilon  \over |r(z,w)|^{1+\epsilon}  r(w)}.
\end{eqnarray*}
 By (2.24):
$$
|\nabla_wY(w,z)|\le {C\over |r(w, z)|^2} \le {C\over |r(z,w)| r(w)},
$$
one has
\begin{eqnarray*}
|{\d \over \d z} \nabla_w H(z,w)|
&\le& Cr(w)\Big( \int_D {1\over |z-\xi|}|\nabla_w X(w,\xi)|dA(\xi)+|\nabla_w Y(w, z)|\Big)\\
&&+C\Big| \int_D {1\over z-\xi} X(w, \xi) dA(\xi)\Big|+|Y(w, z)|\\
&\le &{C_\epsilon \over |r(z,w)|^{1+\epsilon} }.
\end{eqnarray*}
Therefore, (2.16) is proved, and so is the lemma.\epf

We come to study the Green's function for a bounded domain $\Omega$ with $C^{1,\alpha}$ boundary. Choose $z_1,\cdots, z_m\in \d \Omega$ and $\delta>0$ such that 
$\Omega\cap D(z_j, 4 \delta)$ is simply connected and
$$
\d \Omega\subset \cup_{j=1}^n D(z_j, {\delta\over 2} ).
$$
 By the partition of the unity, we choose $\psi_j\in  C^\infty_0(D(z_j, \delta))$ such that
$0\le \psi_j(z)\le 1$,
$$
\psi_0(z)+\sum_{j=1}^m \psi_j(z)=1,\quad z\in \Omegabar
$$
and
$$
\hbox{supp}(\psi_0)\subset \{z\in \Omega: \hbox{dist}(z, \d \Omega)>\delta/10\}.
$$
Choose $ D(z_j,2 \delta)\subset \tilde{\Omega}_j \subset D(z_j, 3\delta)$ with $C^\infty$ boundary such that $\Omega_j=\Omega\cap\tilde{\Omega}_j$
has $C^{1,\alpha}$ boundary and
$$
\Omega\cap D(z_j, 2\delta)\subset \Omega_j\subset \Omega\cap D(z_j, 3\delta),\quad 1\le j\le m. \leqno(2.26)
$$
Let $\phi_j: \Omega_j\to D(0,1)$ be a fixed biholomorphic map. By Proposition \ref{R1}, one has $\phi_j\in C^{1,\alpha}(\overline{\Omega_j})$.
For $\epsilon_0>0$, we let $\Omega_j(2\epsilon_0)$ be the $2\epsilon_0$-neighborhood of $\overline{\Omega_j}$.  We first extend $\phi_j \in
 C^{1,\alpha}(\overline{\Omega_j(2\epsilon_0)}\cap \Omegabar)$ from $\Omegabar_j$. Choose $\eta_j \in C^\infty_0(\Omega_j(2\epsilon_0))$ and
 $\eta_j(z)=1$ on $\Omega_j(\epsilon_0)$. Let
 $$
 \phi_j(w)=\cases{ \phi_j(w) \eta_j(w), \hbox{ if } w\in \Omega_j(2\epsilon_0)\cap \Omega,\cr 0,\qquad\qquad \ \hbox{ if } w\in \Omega\setminus \Omega_j(2\epsilon_0).\cr}\leqno(2.27)
 $$
Then $\phi_j\in C^{1,\alpha}(\Omegabar)$ which is an extension from $\phi_j$ on $\Omegabar_j$.

Using the holomorphic changes of variable by $\phi_j(w)$ and the Poisson kernel for $D(0,1)$, one has
$$
\TP_{\Omega_j} [f]( w) ={1\over 2\pi} \int_{\d \Omega_j} {1-|\phi_j(w)|^2 \over |1-\phi_j(w) \overline{\phi_j(\xi)} |^2} f(\xi) 
|\phi_j'(\xi)| |d\xi|, \  w\in \Omega_j.\leqno(2.28)
$$
Let
$$
H_j(z,w)=\TP_{\Omega_j}[{\psi_j(\cdot)-\psi_j(w)\over \pi} \log|z-\cdot|^2](w),\quad (z, w)\in \Omega\times \Omega_j.\leqno(2.29)
$$
Let
$$
d_j(w)=d_{\Omega_j}(w),\quad d_j(z,w)=\cases{ 1-\phi_j(z)\overline{\phi_j(w)}, \quad z\in \Omega_j, w\in \Omega_j,\cr
z-w,\qquad \qquad z\not\in \Omega_j, w\in \Omega_j.\cr}\leqno(2.30)
$$
Applying Lemma \ref{Lem} and holomorphic map $\phi_j: \Omega_j\to D(0,1)$, one has

\begin{corollary}\label{Cor2}  For $z\in \Omega\times \Omega_j$, $H_j(z,w)$ satisfies
$$
|H_j(z,w)|\le C_\epsilon { d_j(w)  \over |d_j(z,w)|^\epsilon } \log \Big( {1+ d_j (w) \over d_j(w)}\Big), \leqno(2.31)
$$
$$
|\nabla_z H_j(z, w)|\le {C_\epsilon d_j(w) \over |d_j(z,w)|^{1+\epsilon} }\leqno(2.32)
$$
 and
$$
\Big|\nabla_w {\d\over \d z}  H_j(z,w)\Big|\le 
{C_\epsilon \over |d_j(z,w)|^{1+\epsilon} }.\leqno(2.33)
$$
\end{corollary}

For $z\in \Omega$ and $w\in \Omega_j$, we let 
\begin{eqnarray*}
(2.34)\quad  u_j(z,w)&=& \TP_{\Omega_j}\Big[ \psi_j(\cdot)  [{\log |z-\cdot|^2\over \pi}](w)\\
 & =&{1\over 2\pi} \int_{\d \Omega_j} {1-|\phi_j(w)|^2 \over |1-\phi_j(w) \overline{\phi_j(\xi)} |^2}{\psi_j(\xi) \log |z-\xi|^2\over \pi} 
|\phi_j'(\xi)| |d\xi|.
\end{eqnarray*}

Since $D(z_j, 2\delta)\cap \Omega\subset \Omega_j\subset D(z_j, 3\delta)\cap \Omega$, one can see that
$$
|1- \phi_j(w)\overline{\phi_j(\xi)} |\ge {1\over C(\delta)},\quad \xi\in D(z_j, \delta)\cap \Omega \hbox{ and }\ 
w\in \{ \lambda\in \Omega: |\lambda-z_j|>3\delta/2\}.
$$
Since $\phi_j\in C^{1,\alpha}(\Omegabar)$  which is an extension of $\phi_j$ from $\Omegabar_j$, then there is an $0<\epsilon_1< \epsilon_0$ such that
$$
|1- \phi_j(w)\overline{\phi_j(\xi})|\ge {1\over 2 C(\delta)}
$$
for all $ \xi\in D(z_j, \delta)\cap \Omega, w\in \{\lambda\in \Omega\cap \Omega_j(2\epsilon_1)\setminus D(z_j, 4\delta/3)\}.$
Choose $\zeta_j\in C^\infty_0(\Omega_j(2\epsilon_1))$ with $\zeta_j=1$ on $\Omega_j(\epsilon_1)$. We define
$$
u_j(z,w)=\cases{\zeta_j(w) u_j(z,w),\quad \hbox{if } \ \in \Omega_j(2\epsilon_1)\cap\Omega,\cr 0,\qquad\qquad\qquad \hbox{if } w\in \Omega_j(2\epsilon_1)^c\cap \Omega.\cr}
\leqno(2.35)
$$
Then
$$
\|u_j(z,\cdot) \|_{C^{1,\alpha}(\Omegabar \setminus D(z_j, 3\delta/2))} \le C(\delta),\quad z\in \Omega. \leqno(2.36)
$$ 
and
 $$
\|\nabla_z u_j(z,\cdot) \|_{C^{1,\alpha}(\Omegabar \setminus D(z_j, 3\delta/2))}\le C(\delta),\quad z\in \Omega. \leqno(2.37)
$$
Choose $\chi\in C^\infty(\CC)$ such that $\chi(w)=0$ if $w\in D(z_j, 3\delta/2)$ 
and $\chi(w)=1$ if $w\not \in D(z_j,  5 \delta/3)$.  Now since
\begin{eqnarray*}
\Delta'_w u_j(z,w)&=&\Delta'_w (\chi(w) u_j +(1-\chi(w)) u_j)\\
&=&{\d \over \d w} {\d (u_j(z,w)\chi )\over \d \wbar}-\Delta'_w \chi u_j(z,w)-2\Re {\d \chi\over \d w}{\d u_j(z,w)\over \d \wbar}\\
&=&{\d \over \d w} F_j^1(z,w)+F_j^2(z,w)
\end{eqnarray*}
and
$$
\| F^1(z,\cdot)\|_{C^\alpha(\Omegabar)}\le C(\delta)\quad\hbox{and }\ 
\|F^2 (z,\cdot)\|_{L^\infty(\Omega)}\le C(\delta)
$$
uniformly for $z\in \Omega$. Thus,
$$
\cases{\Delta_w (U(z,w)-\sum_{j=1}^m u_j(z, w))
=-{\d \sum_{j=1}^m F^1(z,w) \over \d w} -\sum_{j=1}^m F_j(z,w),\quad w\in \Omega,\cr
 (U(z,w)-\sum_{j=1}^m u_j(z, w))=0,\quad w\in \d \Omega.\cr}
$$
By Theorem 8.34 and (8.90) in \cite{GT}, there is a constant $C_\Omega$ depending only on the diameter of $\Omega$ and $\delta$ such that 
$$
\|(U(z,w)-\sum_{j=1}^m u_j(z, w))\|_{C^{1,\alpha}(\Omegabar)}\le C_\Omega, \quad z\in \Omega
$$
and
$$
\Big\| {\d \over \d z} (U(z,\cdot)-\sum_{j=1}^m u_j(z, \cdot)) \Big\|_{C^{1,\alpha}(\Omegabar)}\le C_\Omega,\quad z\in \Omega.
$$
Let
$$
v(z,w)=-U(z,w)+\sum_{j=1}^m u_j(z, w).\leqno(2.38)
$$
Then $v(z,w), \nabla_z v(z,w)\in C^{1,\alpha}(\Omegabar)$ and
$$
U(z,w)=-v(z,w)+\sum_{j=1}^m u_j(z,w).\leqno(2.39)
$$
Let $G_j$ be the Green's functon of $\Omega_j$. Then we extend the definition of $G_j$ to be defined on  $\Omega \times \Omega_j$ and on $\Omega_j\times \Omega$.
For $(z,w)\in \Omega\times \Omega_j$, since
\begin{eqnarray*}
\lefteqn{u_j(z,w) - \psi_j(w){\log |z-w|^2\over \pi}}\\
&=&\TP_{\Omega_j} \Big[ (\psi_j(\cdot)-\psi_j(w)){\log|z-\cdot|^2\over \pi}\Big](w)
-\psi_j(w) G_j(z, w)\\
&=& H_j(z,w)-\psi_j(w) G_j(z, w),
\end{eqnarray*}
one has
\begin{eqnarray*}
G(z,w)&=& {\log|z-w|^2\over \pi}-U(z, w)\\
&=&v(z,w) +\psi_0(w) {\log|z-w|^2\over \pi} +\sum_{j=1}^m \psi_j(w) {\log |z-w|^2\over \pi} -u_j(z,w))\\
&=&v(z,w)+\psi_0(w) {\log|z-w|^2\over \pi}  +\sum_{j=1}^m \psi_j(w) G_j(z,w) -\sum_{j=1}^m H_j(z,w)).
\end{eqnarray*}

\noindent {\bf Note.}  Let
$$
{\cal H}(z, w)=:\sum_{j=1}^m H_j(z,w).\leqno(2.40)
$$
Our  $H_j(z,w)$ was defined on $\Omega\times \Omega_j$, now we have extended it to be defined on $\Omega\times \Omega$ 
through the extension of $\phi_j$ from $\Omega_j$ to $\Omega$ given by (2.27) and the extension of $u_j(z,w)$ from $\Omega_j$ to $\Omega$
given by (2.35). Combining these extensions and Corollary 2.4, we have the following estimates.
$$
|{\cal H}(z,w)|\le C_\epsilon {d_\Omega(w) \over d(z, w)}| \log d_{\Omega}(w)|,\leqno (2.41)
$$
$$
|\nabla_z {\cal H}(z,w)|\le C_\epsilon {d_\Omega(w)\over d(z,w)^\epsilon}\leqno(2.42)
$$
and
$$
\Big|\nabla_w {\d {\cal H}\over \d z}\Big|\le {C_\epsilon \over d(z,w)^{1+\epsilon}},\leqno(2.43)
$$
where
$$
d(z, w)=d_\Omega (w)+|z-w|.\leqno(2.44)
$$
Therefore, we have proved the following theorem.

\begin{theorem}\label{Green} Let $\Omega$ be a bounded domain in $\CC$ with $C^{1,\alpha} $ boundary for some $\alpha>0$.
 Let $G(z,w)$ be the Green's function for $\Delta$ in $\Omega$.
Then
$$
(2.45)\  \ G(z,w)=\psi_0(w){ \log |z-w|^2\over \pi}
+\sum_{j=1}^m \psi_j(w) G_j(z,w)
-{\cal H}(z,w)+v(z,w),
$$
where $v(z, w),  \nabla_z v(z,w)\in C^{1,\alpha}( \Omegabar)$ as functions of $w$ and $H$  is given by (4.40) satisfy the estimates (2.41)--(2.44).
\end{theorem}

\begin{corollary} Let $\Omega$ be a bounded domain in $\CC$ with $C^{1,\alpha}$ boundary for some $\alpha>0$. Then
the Bergman projection ${\cal P}$ is bounded on $L^p(\Omega)$ for $1<p<\infty$.
\end{corollary}

\proof Since
$$
K(z,w)={\d^2 G(z,w)\over \d z \d\wbar},\quad z\ne w,\leqno(2.46)
$$
Let
$$
V(z,w)={1\over \pi (z-w)} {\d \psi_0(w)\over \d \wbar}+{\d^2 v(z,w)\over \d z\d\wbar}.\leqno(2.47)
$$
By Theorem \ref{Green}, we have
\begin{eqnarray*}
\lefteqn{\ (2.48)\quad K(z,w)}\\
&\quad=&V(z,w)+\sum_{j=1}^m {\d^2 \psi_j(w) G_j(z,w)\over \d z\d\wbar}-{\d^2 {\cal H}(z,w) \over \d z\d \wbar}\\
&\quad=&V(z,w)+\sum_{j=1}^m {\d\psi_j(w) \over \d \wbar} {\d G_j(z,w)\over \d z }
+\sum_{j=1}^m \psi_j(w){\d^2 G_j(z,w)\over \d z\d\wbar} -{\d^2 {\cal H} (z,w) \over \d z\d \wbar}.
\end{eqnarray*}
We know that 
$$
|V(z,w)+\sum_{j=1}^m {\d\psi_j(w) \over \d \wbar} {\d G_j(z,w)\over \d z }|\le {C\over |w-z|}.
\leqno(2.49)
$$
By (2.33) with $\epsilon=1/2$ and Green's function $G_j$,  one has
$$
 |\sum_{j=1}^m \psi_j(w){\d^2 G_j(z,w)\over \d z\d\wbar} -{\d^2 {\cal H} (z,w) \over \d z\d \wbar}|\
\le {C\over |d(z,w)|^2},\leqno(2.50)
$$
where $d(z,w)$ are given by (2.44). For any $1<p<\infty$, choose $\epsilon>0$ such that $p\epsilon, p' \epsilon\le 1/2$. 
One can easily see that
$$
\int_\Omega |K(z,w)| r(w)^{-p\epsilon} dA(w)\le {C\over p\epsilon} d_\Omega(z)^{-p\epsilon}
$$
and
$$
\int_\Omega |K(z,w)| r(w)^{-p'\epsilon} dA(w)\le {C\over p'\epsilon} d_\Omega(z)^{-p'\epsilon}
$$
By the Schur's lemma, we have the Bergman projection ${\cal P}[f]=\int_\Omega K(z,w) f(w) dA(w)$ is bounded on the both $L^p$ and
$L^{p'}$ with norm less than or equal to ${C\over \epsilon}$. Therefore, the corollary is proved.\epf

\section{Formula solution to $\dbar$-equations}

Let $\rho\in C^{1,\alpha}(\CC)$ and $\Omega=\{z\in \CC: \rho(z)<0\}$. There is a scalar constant $C_0>1$ such that
${1\over C_0}\le |\nabla \rho(z)|\le C_0$ on $\d \Omega$. Let $C_\Omega$ denote a constant depending only on $C_0$
and $\|\rho\|_{C^{1,\alpha}(\Omegabar)}$. It may not be the same at each appearance.

Let $G=G_\Omega$ be the Green's function for $\Delta'={\d^2 \over \d z\d \zbar}$ on $\Omega$. Define
$$
k(z, w)={\d G_{\Omega}(z, w)\over \d z}
 \leqno(3.1)
$$
and
$$ \label{S1}
T[f](z)=\int_\Omega k(z,w) f(w) dA(w).\leqno(3.2)
$$


\begin{proposition} Let $\Omega $ be a bounded domain in $\CC$ with $C^{1,\alpha}$ boundary for some $\alpha>0$ and $1\le p\le \infty$. Then

(i)   $T[f]  $ is 
 the canonical solution of  $\dbar u=f d \zbar$;
 
 (ii) $T: L^p(\Omega)\to L^{2p \over 2-p}(\Omega)$ if $1\le  p<2$;

 (iii)  $T: L^p(\Omega)\to L^\infty(\Omega)$ is bounded if $2<p\le \infty$;
 
 (iv) $T: L^p(\Omega)\to C^{1-2/p}(K)$  for any compact set $K\subset 
 \Omega$.
\end{proposition}

\proof By (2.2), the definition of $T[f]$ and the definition of the Green's function, one can easily see that 
$$
{\d T [f]\over \d \lambdabar} (\lambda)={\d^2 G[f] \over \d \lambda\d \lambdabar} =f(\lambda),\quad \lambda \in \Omega.
$$
For any $z_0\in \Omega$,
since $G(z, w)=0$ on $\d \Omega\times \Omega$,  by Proposition \ref{R1}, Theorem \ref{Green},   and  integration by parts, one has
$K(z_0, \cdot)\in C^\alpha(\Omegabar)$ and
$$
\int_{\d \Omega} G(\lambda, w) K(z_0, \lambda) d\lambda=0,\quad w\in \Omega.
$$
Thus, by (3.1) and $K(z_0, \cdot)\in C^\alpha(\Omegabar)$, for $w\in \Omega$, one has
$$
\int_\Omega k(\lambda, w) K(z_0, \lambda) dA(\lambda)
= \int_{\d \Omega} G(\lambda, w) K(z_0, \lambda) {i d\lambdabar \over 2} - \int_\Omega G(\lambda, w) {\d K(z_0, \lambda) \over \d \lambda} dA(\lambda) =0.
$$
By the above identity, one has 
\begin{eqnarray*}
\int_\Omega T[f](\lambda) K(z_0, \lambda)  dA(\lambda)
&=&\int_\Omega \int_\Omega k(\lambda, w) K(z_0, \lambda) dA(\lambda) f(w) dA(w)\\
&=&\int_\Omega 0\cdot  f(w) dA(w)\\
&=& 0.
\end{eqnarray*}
Therefore, $T[f]\in A^2(\Omega)^{\perp}$ sinc $K$ is the Bergman kernel of $\Omega$. So, $T[f]$ is the canonical solution of $\dbar u =f d\zbar $ in $\Omega$. Part (i) is proved.

\medskip

Let
$$
v(z)={1\over \pi}\int_\Omega {f(w)\over z-w} dA(w).\leqno(3.3)
$$
Then ${\d v\over \d \zbar}=f$. If $f\in L^p(\Omega)$ for $1\le  p<2$, we have $v\in W^{1,p}(\Omega)$. By the Sobolev embedding theorem,
one has $v\in L^{2p\over 2-p}(\Omega)$. Since the Bergman projection ${\cal P}$ is bounded on $L^q(\Omega)$ for any $1<q<\infty$ by
Corollary 2.6. This implies that $T[f]=v-{\cal P}[v]\in L^{2p\over 2-p}(\Omega)$. Therefore, Part (ii) is proved.

\medskip

By Proposition \ref{R1} and Theorem \ref{Green},  one has
 $$
 |k(z, w)|=|{\d G(z, w) \over \d z}| \le {C_\Omega \over |z-w|}.\leqno(3.4)
 $$
 This implies that if $2<p\le \infty$, then
 $$
 |T[f](z)|\le C_\Omega\int_\Omega{ |f(w)|\over |w-z|} dA(w)\le {C_\Omega\over 2-p'} \|f\|_{L^p} \le C_\Omega {p-1\over p-2}  \|f\|_{L^p}.
 $$
 This means $\|T[f]\|_{L^\infty}\le C_\Omega {p-1\over p-2} \|f\|_{L^p}$ if $p>2$. Part (iii) is proved.
\medskip

 By Sobolev embedding theorem, one has that  $v\in W^{1, p}(\Omega)\subset C^{1-2/p}(\Omegabar)$ for $2< p<\infty$.  Thus,
$$
T[f]=v-{\cal P}[v]\in C^{1-2/p}(K),\quad\hbox{for any compact set } K\subset \Omega.
$$
So, the proof of Part (iv) is completed.
Therefore, the proof of the proposition is complete.\epf
\medskip

\noindent  {\bf Remark.}  For (iv), one expects that
 $T: L^p(\Omega)\to C^{\alpha_0}(\Omegabar)$, where $\alpha_0=\min\{\alpha, 1-2/p\}$.
 This can be done easily when $\Omega$ is simply connected. For a general case, one needs to do some more works from
 estimations of Green's function given by Theorem \ref{Green}. Since we don't use this property, we omit here.
\medskip

For any  $1\le j\le n$ and $z\in \CC^n$, write
$$
z^{(j)}=(z_1,\cdots, z_{j-1}, z_{j+1}, \cdots, z_n), \quad z=(z_j; z^{(j)}).\leqno(3.5)
$$
Let $f\in L^2(\Omega^n)$, we define
the Bergman projection $P_j: L^2(\Omega)\to A^2(\Omega)$  by
$$
P_j f (z)={\cal P} [f(\cdot, z^{(j)})](z_j)
=\int_\Omega K(z_j, w_j) f(w_j; z^{(j)}) dA(w_j), \leqno(3.6)
$$
for almost every $z^{(j)}\in \Omega^{n-1}.$ We also use the notations $P_0=P_{n+1}=I$. Similarly,  we also use the following notation:
$$
T_j f(z)=T[f(\cdot\, ; z^{(j)})] (z_j), \quad 1\le j\le n.\leqno(3.7)
$$

The following theorem is a very important formulation for the canonical solution of $\dbar u=f$.

\begin{theorem} Let $\Omega$ be a bounded domain in $\CC$ with $C^{1,\alpha}$ boundary for some $\alpha>0$. For $1\le  p\le \infty$ and any $\dbar$-closed $(0,1)$-form
 $f=\sum_{j=1}^n f_j d\zbar_j \in L^p_{(0,1)}(\Omega^n)$, the canonical solution 
$u=S[f] \in L^p (\Omega^n)$  to $\dbar u=f$ satisfies
$$
S[f](z)=\sum_{j=1}^n T_j  P_{j-1}\cdots P_{0} f_j=\sum_{j=1}^n T_j P_{j+1}\cdots P_{n+1} f_j. \leqno(3.8)
$$
\end{theorem}

\noindent  {\bf Note.} $S[f] \in L^p(\Omega^n)$ is the canonical solution here means $P_{\Omega^n} [S[f]]=0$.
\medskip

\proof  For each $1\le j\le n$, since ${\d u(z_j; z^{(j)})\over \d \zbar_j} =f_j ( z_j; z^{(j)})\in L^p(\Omega)$,
by the estimates on the Green's function given  by  Proposition 2.2 and Theorem \ref{Green}, one has that
$$
u(z_j; z^{(j)})-P_j[u(\cdot\, ; z^{(j)})](z_j)=T_j[f_j(\cdot\, ; z^{(j)})](z_j),\leqno(3.9)
$$
for almost every $ z^{(j)}\in \Omega^{n-1}.$ 

  Since $u-P_1[u]$ is the canonical solution of ${\d u\over \d \zbar_1}=f_1$, one has
$$
P_0 u-P_1 P_0 u=u-P_1[u]=T_1 f_1=T_1 P_0 f_1.
$$
Similarly, 
$$
P_1P_0[u]-P_2P_1P_0[u]=P_1[(I-P_2) u]=P_1T_2{[\d u\over \d\zbar_2}] = P_1T_2[ f_2]=T_2P_1[f_2].
$$ Keeping the same process, one has
$$
P_{j-1}\cdots P_1 P_0u-P_j P_{j-1} \cdots P_1P_0 u= T_j P_{j-1}\cdots P_1  f_j,\quad 1\le j\le n.
$$
 Since $P_1\cdots P_n u=P_{\Omega^n} u=0$ and $P_0=I$, one has
$$
S[f]=u=\sum_{j=1}^n (P_{j-1}\cdots P_0u-P_j P_{j-1} \cdots P_0 u )=\sum_{j=1}^n  T_j P_{j-1}\cdots P_{1} P_0 f_j.
$$

On the other hands, let $P_{n+1}=I$, then
$$
 u-P_n  u=T_n f_n.
$$
With the same process, one has
$$
P_n \cdots P_j u-P_n\cdots P_j  P_{j-1}  u=T_{j-1} P_j\cdots P_n f_{j-1}.
$$
Since $u$ is the canonical solution of $\dbar u=f$, one has $P_{n+1} P_n\cdots P_1 u=0$ and
$$
\sum_{j=1}^{n} T_{j} P_{j+1}\cdots P_{n+1} f_{j}=\sum_{j=1}^{n} (P_{n+1} P_n \cdots P_{j+1} u- P_{n+1} P_n\cdots P_j   u)=u.
$$
These prove (3.8), so, the proof of  Theorem 3.2  is complete. \epf

\begin{theorem}  Let $1<p<\infty$ and  let $\Omega$ be a bounded  domain in $\CC$ with $C^{1,\alpha}$ boundary for some $\alpha>0$.
If   $f_m , f\in L^p_{(0,1)}(\Omega^n)$ with $f_m\to f$ in $L^p_{(0,1)}(\Omega^n)$, then
$$
\lim_{m\to \infty}\Big \|S[f_m]-S[f]\Big\|_{L^p(\Omega^n)}=0.\leqno(3.10)
$$
\end{theorem}

\proof Notice that $T_j$ and $P_j$ are bounded on $L^p(\Omega)$ with finite norms $\|T_j\|_p$ and $\|P_j\|_p$. 
Write
$$
f_m=\sum_{k=1}^n (f_m)_k d\zbar_k,\quad f=\sum_{k=1}^n f_k d\zbar_k.
$$
By (3.8), one has
\begin{eqnarray*}
\|S[f_m]-S[f]\|_{L^p(\Omega^n)}&\le& \sum_{j=1}^n\|T_j P_{j-1}\cdots P_0[(f_m)_j-f_j]\|_{L^p(\Omega^n)}\\
&\le&\sum_{j=1}^n \|T_j\|_p \|P_{j-1}\|_p\cdots \|P_1\|_p \|(f_m)_j-f_j\|_{L^p(\Omega^n)}\\
&\le&\sum_{j=1}^n \|T_j\|_p \|P_{j-1}\|_p\cdots \|P_1\|_p \|f_m-f\|_{Lp_{(0,1)}(\Omega^n)}
\end{eqnarray*}
which converges to $0$ as $m\to\infty$. The proof is complete.\epf

\section{Regularity and  a new formula solution}

To prove a uniform $L^p$-estimate for $S[f]$ when $p\ge 1$. By viewing 
$$
T_jP_{j-1}[f]={1\over \pi} \int_\Omega \int_\Omega k(z_j,w_j) K(z_{j-1},w_{j-1}) f(w) dA(w),
$$
one see that the singularity of $k(z_j,w_j)$ is about ${1\over |w_j-z_j|}$ and $|K(z_{j-1}, w_{j-1}) |$ is about $ {1\over d_\Omega(w_{j-1})^2+|z_{j-1}-w_{j-1}|^2}$.
In this section, through the integration by part, we try to move the higher singularity partially to the lower singularity. In the end, we can average
them to get kernel on $\Omega\times \Omega$ with the singularity like ${1\over |z_j-w_j|^{3/2}|z_{j-1}-w_{j-1}|^{3/2}}$. Which becomes integrable on $\Omega^2$.
For this purpose, we need to introduce some notation and to do preparations. First, we state the following lemma.

\begin{lemma} Let $(X, d \mu)$ be a measurable space with positive Borel measure $d\mu$. 
Let ${\cal  K}$ be measurable function on $X \times X$ such that
$$
\int_{X} |{\cal K}(z, w)| d\mu (w) +\int_{X} |{\cal K}(z, w)| d\mu (z)\le C, \quad\hbox{ for all }  z, w\in X.\leqno(4.1)
$$
Then the integral operator
$$
{\cal T} [f](z)=\int_X {\cal K}(z,w) f(w) d\mu (w)\leqno(4.2)
$$
is bounded on $L^p(X)$ and $\|T\|\le C$ for all $1\le p\le \infty$.
\end{lemma}
\proof The proof can be followed from Schur's lemma with test function $1$. One can also prove it by using the interpolation theorem.
Since ${\cal T}$ is bounded on $L^1$ with the norm $\le C$ and it is also bounded on $L^\infty$ with norm $\le C$. The interpolation theorem implies
that ${\cal T}$ is bounded on $L^p$ with the norm $\le C$.\epf


\medskip

For any $1\le i \ne j\le n$, define
$$
\tau_{i,j}=|z_i-w_i|^2+|z_j-w_j|^2 \quad \hbox{and }\quad \tau^k_{i,j}={|w_k-z_k|^2\over |w_i-z_i|^2+|w_j-z_j|^2}.\leqno(4.3)
$$
Define $\d_{\jbar}={\d \over \d \wbar_j}$ and
\begin{eqnarray*}
\lefteqn{ (4.4)\quad b^{i, j}(z,w) }\\
&\qquad :=&\d_{\jbar} \Big(k(z_j, w_j) \tau^j_{i,j}\Big)\\
&\qquad =& k(z_j , w_j) \d_{\jbar} \tau_{i,j}^j +\tau^j_{i,j} \d _{\jbar} k(z_j, w_j)\\
&\qquad =& k(z_j, w_j) {(w_j-z_j) |w_i-z_i|^2\over (\tau_{i,j})^2 }+ \tau_{i,j}^j  \d_{\jbar}  k(z_j, w_j)\\
&\qquad =& k(z_j, w_j) {(w_j-z_j) \over \tau_{i,j} }- {k(z_j, w_j) (w_j-z_j) \over \tau_{i,j } } \tau^j_{i,j}+{|w_j-z_j|^2\over \tau_{i,j}}\d_{\jbar} k(z_j, w_j)\\
&\qquad =&\Big(h(z_j, w_j)+(\wbar_j-\zbar_j) (\d_{\jbar} h(z_j, (w_j)\Big) {1 \over \tau_{i,j}}-{h(z_j, w_j)\over \tau_{i,j}} \tau^j_{i,j},
\end{eqnarray*}
where
$$
h(z_j, w_j)=(w_j-z_j)k(z_j, w_j)
.\leqno(4.5)
$$
By Proposition 2.2 and Theorem 2.5, one has
$$
|k(z_j, w_j)|\le {C_\Omega \over |z_j-w_j|}\leqno(4.6)
$$
and
$$
 |k(z_j, w_j)|\le C_\Omega { d_\Omega(w_j)\over |z_j-w_j|^2} .\leqno(4.7)
$$
These imply that
$$
|h(z_j, w_j)+ |(\wbar_j-\zbar_j) {\d h(z_j, w_j)\over \d \wbar_j}| \le C_\Omega  \ \hbox{ and } \  |h(z_j, w_j)|\le {C_\Omega d_\Omega(w_j) \over |z_j-w_j|}. \leqno(4.8)
$$
Therefore
$$
|b^{i,j}(z, w)|\le {C_\Omega\over \tau_{i,j}(z,w)}.\leqno(4.9)
$$
Thus,
\begin{eqnarray*}
\lefteqn{(4.10)\qquad |k(z_j, w_j)||b^{j,i}(z,w)|+|k(z_i, w_i)||b^{i,j}(z,w)| }\\
&\qquad \qquad\le & {C_\Omega \over |z_j-w_j| \tau_{i,j}}+{C_\Omega \over |w_i-z_i|\tau_{i,j} }  \qquad\qquad\qquad\qquad\qquad\qquad\qquad\qquad\qquad \\
&\qquad \qquad \le & {C_\Omega  \over |w_i-z_i|^{3/2} |w_j-z_j|^{3/2}}.
\end{eqnarray*}
By (4.7) and integration by parts, one has
\begin{eqnarray*}
T_j P_i  [f_j ] &=&T_j[f_j]-T_j T_i [{\d f_j \over \d \wbar_i}]\\
&=&T_j[f_j] +\int_{\Omega^2} k_j {\d (k_i \tau^i_{i,j} ) \over \d \wbar_i} f_j
+\int_{\Omega^2} k_i {\d (k_j \tau^j_{i,j} ) \over \d \wbar_j} f_i\\
&=&T_j[f_j] +\int_{\Omega^2} k_j  b^{j, i} f_j
+\int_{\Omega^2} k_i b^{i,j} f_i.
\end{eqnarray*}
By Lemma 4.1, estimations (4.6), (4.10) and the above identity, we have proved the following proposition.

\begin{proposition} Let $\Omega$ be a bounded domain with $C^{1,\alpha}$ boundary. 
Let $f\in C^1_{(0,1)}(\Omegabar^n)$ be $\dbar$-closed. Then
$$
\|T_j P_i [f_j]\|_{L^p(\Omega^n)}
\le C_\Omega (\|f_j\|_{L^p(\Omega^n)}+(\|f_i\|_{L^p(\Omega^n)}\Big), \leqno(4.11)
$$
for all $1\le p<\infty$, where $C_\Omega$ is a constant defined in the beginning of Section 3.
\end{proposition}

Let $k_j=k(z_j, w_j)$. Then
\begin{eqnarray*}
T_j P_\ell P_i[ f_j]
&=&T_j P_\ell [f_j] +P_\ell \Big[\int_{\Omega^2} k_j  b^{j, i} f_j
+\int_{\Omega^2} k_i b^{i,j} f_i\Big]\\
&=&T_j P_\ell [f_j] +\int_{\Omega^2} k_j  b^{j, i} f_j
+\int_{\Omega^2} k_i b^{i,j} f_i\\
&&- T_\ell \Big[\int_{\Omega^2} k_j  b^{j, i} {\d f_j\over \d\wbar_\ell}\Big]-
T_\ell \Big[\int_{\Omega^2} k_i b^{i,j}{\d  f_i \over \d \wbar_\ell}\Big]\\
&=&T_j P_\ell [f_j] +T_j P_i [f_j] -T_j[f_j] \\
&&- T_\ell \Big[\int_{\Omega^2} k_j  b^{j, i} {\d f_j\over \d\wbar_\ell}\Big]-
T_\ell \Big[\int_{\Omega^2} k_i b^{i,j}{\d  f_i \over \d \wbar_\ell}\Big].
\end{eqnarray*}
The estimate (4.11) holds for the first and the second in the above identity. Their term is similar to the fourth term. We need only to consider
$T_\ell \Big[\int_{\Omega^2} k_j  b^{j, i} {\d f_j\over \d\wbar_\ell}\Big]$.
 Since $\ell\ne i, j$, one has
 $$
 \d_{\lbar} \Big( k_\ell \tau^\ell_{j,\ell} k_j  b^{j,i} )
 =k_j b^{j,\ell} b^{j, i},\leqno(4.12)
 $$
 and
 $$
 \d_{\jbar} \Big( \tau^j_{j,\ell} k_\ell k_j  b^{j,i})=k_\ell b^{\ell, j} b^{j,i} -k_\ell {1\over \tau_{\ell, j} }a^{j, i},\leqno(4.13)
 $$
 where
$$
 a^{i, j}=|z_j-w_j|^2 k_j \d_{\jbar}( b^{j, i}).\leqno(4.14)
 $$
  Since
$$
(4.15)\ \d_{\ibar}b^{i,j}
=\Big(  -b^{i,j} +{h(z_, w_j)\over \tau_{i,j}} \tau_{i,j}^j  \Big) \d_{\ibar}  \log \tau_{i,j}
=\Big(  -b^{i,j} +{h(z_, w_j)\over \tau_{i,j}} \tau_{i,j}^j  \Big) {(w_i-z_i)\over  \tau_{i,j}}.
$$
By (4.6), one has
$$
\Big|{\d b^{i,j}\over \d \wbar_i}\Big|
\le {|z_i-w_i|\over ( \tau_{i,j} )^2}.\leqno(4.16)
$$
Thus,
$$
|a^{i,j}|=| k_j (z_j-w_j)^2 \d_{\jbar} b^{j,i}|
\le |k_j |z_j-w_j|^2 { |z_j-w_j|\over |\tau_{i,j}|^2}
\le {C_\Omega\over \tau_{i,j}}.\leqno(4.17)
$$
\begin{eqnarray*}
(4.18)\  T_\ell [\int_{\Omega^2} k_j b^{j,i}{\d f_j\over \d \wbar_\ell} ]
&=&\int_{\Omega^3} k_\ell  k_j b^{j,i}{\d f_j\over \d \wbar_\ell} \\
&=&-\int_{\Omega^3} \d_{\lbar} \Big( k_\ell \tau^\ell_{j,\ell} k_j  b^{j,i} \Big) f_j
-\int_{\Omega^3} \d_{\jbar} \Big( \tau^j_{j,\ell} k_\ell k_j  b^{j,i} \Big) f_\ell\\
&=&-\int_{\Omega^3} (k_j  b^{j,\ell} b^{j,i} f_j -k_\ell  b^{\ell, j} b^{j, i} f_\ell )-\int_{\Omega^3} { k_\ell  \over \tau_{j,\ell} } a^{j,i}   f_\ell. \qquad\qquad\qquad
\end{eqnarray*}
Let
$$
E^1_{j, i}=b^{j,i}, \quad E^2_{j, i} =b^{j,i} \ \hbox{ or }\ a^{j, i}\leqno(4.19)
$$
Then
$$
|E^2_{j, i}|\le  {C_\Omega \over |z_j-w_j|^{ 1+\epsilon}  |w_i -z_i|^{2-\epsilon}}.\leqno(4.20)
$$
Using the inequality
$$
a^{2/p} b^{2/q} \le {1\over p} a^2 +{1\over p'} b^2\le (a^2+b^2),\leqno(4.21)
$$
with $n\epsilon=1$,  one has
$$
|k_j b^{j,\ell} \,  E^1_{j,i}| \le {C_\Omega\over |z_j-w_j|}  {1\over |z_j-w_j|^{2\epsilon}} {1\over |z_i-w_i|^{2-\epsilon} |w_\ell-z_\ell|^{2-\epsilon}}\leqno(4.22)
$$
and
$$
|k_\ell {1\over \tau_{j,\ell}}E^2_{j, i}|\le { C_\Omega \over |w_\ell-z_\ell|} {1\over |z_\ell-w_\ell|^{2\epsilon} }{1\over |z_j-w_j|^{2-\epsilon} |w_i-z_i|^{2-\epsilon}}.\leqno(4.23)
$$
By Lemma 4.1 and the above estimates, one has proved
$$
\|T_jP_\ell P_i[f_j]\|_{L^p(\Omega^n)}\le C_\Omega \|f\|_{L^p_{(0,1)}(\Omega^n)}.\leqno(4.24)
$$
Let
$$
c^{j, \ell} =k_\ell |z_\ell-w_\ell |^2 {\d \over \d \wbar_\ell} \Big( {1\over \tau_{\ell, j}}\Big)=k_\ell |w_\ell-z_\ell |^2 {z_\ell- w_\ell \over (\tau_{\ell, j})^2}.\leqno(4.25)
$$
Then
$$
|c^{j,\ell} |\le {C_\Omega\over \tau_{j,\ell}}.\leqno(4.26)
$$
For $I=\{i_1,\cdots, i_k\}$ consisting of $k$ different numbers in $\{1,2,\cdots, n\}$, we define
$E_{j, I}$
to the $k$-pairs products from elements in $\{E^1_{i,j} , E^2_{i,j}, c^{i,j}, {1\over \tau_{i,j}}\}$ such that
and 
$$
|E_{j, I}|\le {C_\Omega \over |z_j-w_j|^{k\epsilon} }{1\over \prod_{\ell=1}^k |z_\ell-w_\ell|^{2-\epsilon}}.\leqno(4.27)
$$
Then
$$
P_q[\int_{\Omega^{k+2} } k_j b^{j, \ell}  E_{\ell, I} f_j ]
=\int_{\Omega^{k+2} } k_j b^{j, \ell} E^2_{j, I} f_j -T_q [\int_{\Omega^{k+2} } k_j b^{j,\ell} E_{\ell, I} \d_{\qbar} f_j ]
$$
The function $F$ defined by the  first term has been proved in $L^p$ with $\|F\|_{L^p}\le C_\Omega \|f\|_{L^p_{(0,1)}}$. The second term
\begin{eqnarray*}
\lefteqn{-T_q [\int_{\Omega^{k+2} } k_j b^{j,\ell} E_{\ell, I} \d_{\qbar} f_j ] } \\
&=&\int_{\Omega^{3+k} } k_j b^{q, j} b^{j, \ell} E_{\ell, I} f_j
+\int_{\Omega^{k+3} }k_q  b^{j, q} b^{j, \ell} E_{\ell, I} f_q+\int_{\Omega^{3+k} }k_q {1\over  \tau_{q, j} }a^{j,\ell} E^2_{\ell, I} f_q \\
&=&\int_{\Omega^{3+k}}k_j b^{q,j} E_{j, J} f_j
+\int_{\Omega^{3+k}}k_q b^{j,q} E_{j, J} f_q+\int_{\Omega^{3+k}}k_q {1\over \tau _{q, j}} E_{j, J} f_q,
\end{eqnarray*}
where $J=\{j, I\}$ and $|J|=k+1$. It is easy to see that $E_{j, J}$ satisfies (4.27) by replacing $k$ by $k+1$.
The function $F$ defined by the above three integrals  are in $L^p(\Omega^n)$ with 
$\|F\|_{L^p(\Omega^n}\le C_\Omega \|f\|_{L^p_{(0,1)}(\Omega^n)}$ by applying Lemma 4.1 and (4.21).
Moreover, 
\begin{eqnarray*} 
\lefteqn{-T_q[\int_{\Omega^{k+2}} k_j {1\over \tau_{j,\ell}} E_{\ell, I} \d_{\qbar} f_j]}\\
&=&\int_{\Omega^{3+k}}k_j b^{j, q} {1\over \tau_{j, \ell}} E_{\ell, I} f_j
+\int_{\Omega^{3+k}} k_q b^{q, j} {1\over \tau_{j, \ell}} E_{\ell, I} f_q
+\int_{\Omega^{3+k}}k_j k_q \tau^j_{j,q} \d_{\jbar} {1\over \tau_{j,\ell}} E_{\ell, I} f_q\\
&=&\int_{\Omega^{3+k}}k_j b^{j, q} E_{j, J} f_j +\int_{\Omega^{3+k}}k_q b^{q,j} E_{j, J} f_q+\int_{\Omega^{3+k}}
k_q {1\over \tau_{j, q}}E_{j, J} f_q,
\end{eqnarray*}
where $E_{j, J}$ staisfies (4.27) by replacing $k$ by $k+1$. Applying Lemma 4.1 and  (4.21), one can show
that the function $F$ defined by the above three integrals satisfy
the estimate $\|F\|_{L^p(\Omega^n)}\le C_\Omega\|f\|_{L^p_{(0,1)}(\Omega^n)}$ by applying Lemma 4.1. 

As a summary,   we have proved the following theorem.

\begin{theorem} Let $\Omega$ be a bounded domain in $\CC$ with $C^{1,\alpha}$ boundary for some $\alpha>0$.
Then there is a constant $C=C_{\Omega}$ defined in the beginning of Section 3
 such that for any $\dbar$-closed $(0,1)$-form  $f\in C^1_{(0,1)}(\Omegabar^n)$, one has 
$$
\|S[f](z)\|_{L^p(\Omega)}\le C_\Omega \|f\|_{L^p_{(0,1)}(\Omega)},\quad \hbox{for all }1\le p\le \infty. \leqno(4.28)
$$
\end{theorem}

\section{Proof of Theorem 1.1} 

\subsection{Approximation}

\begin{theorem} Let $\Omega$ be a bounded simply connected domain in $\CC$ with $C^{1,\alpha}$ boundary for some $\alpha>0$. For any $1\le p<\infty$ and
$f\in L^p_{(0,1)}(\Omega^n)$ be $\dbar$-closed, then there is a  $\dbar$-closed squence $\{f_m\}_{m=1}^\infty \subset C^1_{(0,1)}(\Omegabar^n)$
such that 
$$
\lim_{m \to\infty} \|f_m-f\|_{L^p_{(0,1)}}=0.\leqno(5.1)
$$
\end{theorem}

\proof When $\Omega$ is the unit disk $D$, let $\chi^j\in C^\infty_0(D)$ be nonnegative and $\int_{D}\chi^j dA=1$. 
Let $\chi^j_\epsilon=\chi^j(z/\epsilon) \epsilon^{-2}$ and $\chi_\epsilon(z)=\chi^1_\epsilon \cdots \chi^n_\epsilon$
on $D^n$. The proof for this case  is very simple. For any $0<r<1$ and
$\epsilon=(1-r)/2$, since $f_r(z)=f(rz)$ is $\dbar$-closed in $D(0, 1/r)$ and then
$$
F_r(z)=f_r* \chi_\epsilon \in C^\infty_{(0,1)}(\overline{D}^n) \leqno(5.2)
$$
is $\dbar$-closed in $D^n$ and
$$
\|F_r-f\|_{L^p_{(0,1)}(D^n)}\to 0 \leqno(5.3)
$$
as $r\to 1^-$ and any $p\in [1,\infty)$. This argument remains true when $\Omega$ is a simply connected domain
in $\CC$ with $C^{1,\alpha}$ boundary for any $0<\alpha<1$. Let $\phi: \Omega\to D$ be a biholomorphic mapping. Then
$\phi\in C^{1,\alpha}(\Omegabar)$, and $\Omega=\phi^{-1}(D)$, with slightly modification of the unit disc case, one can
similarly prove the theorem.
 \epf

\medskip

Now we are ready to prove Theorem 1.1 when $\Omega$ is bounded simply connected with $C^{1,\alpha}$ boundary. 

\subsection{Proof of Theorem 1.1 if $\Omega$ is simply connected}

\proof For any $1<p<\infty$, by Theorem 4.1, there is a sequence $\{f_m \}_{m=1}^\infty\subset C^1_{(0,1)}(\Omegabar)$ which are
$\dbar$-closed such that
$$
\lim_{m\to\infty} \|f_m-f\|_{L^p_{(0,1)}(\Omega)}=0.\leqno(5.4)
$$
By estimations obtained in Section 4, one has that
$$
\dbar S[f_m]=f_m \leqno(5.5)
$$
and $S[f_m]$ is a canonical solution. Moreover,
$$
\lim_{m\to\infty} \|S[f_m]-S[f]\|_{L^p(\Omega^n)}=0. \leqno(5.6)
$$
For $1< p<\infty$, by Corollary 2.6, one has
\begin{eqnarray*}
\lefteqn{\|S[f]\|_{L^p(\Omega^n)}}\\
&\le & \|S[f_m]\|_{L^p (\Omega^n) } + \|S[f_m]-S[f]\|_{L^p (\Omega^n)}\\
&\le& C_\Omega  \|f_m\|_{L^p_{(0,1)}(\Omega^n)} + \|S[f_m]-S[f]\|_{L^p (\Omega^n)}\\
&\le & C_\Omega \|f\|_{L^p_{(0,1)}(\Omega^n)} +C_\Omega \|f_m-f\|_{L^p_{(0,1)}(\Omega^n)} + \|S[f_m]-S[f]\|_{L^p (\Omega^n)},
\end{eqnarray*}
where $C_\Omega$ is a constant depends neither on  $m$ nor $p$.  Let $m\to \infty$, one has
$$
\|S[f]\|_{L^p_{(0,1)}(\Omega^n)}\le C_\Omega\|f\|_{L^p_{(0,1)}(\Omega^n)}, \quad 1< p<\infty. \leqno(5.7)
$$
Letting $p\to +\infty$, one has
$$
\|S[f]\|_{L^\infty_{(0,1)}(\Omega^n)}\le C_\Omega \|f\|_{L^\infty_{(0,1)}(\Omega^n)}.\leqno(5.8)
$$
The proof of Theorem 1.1 is complete if $\Omega$ is simply connected with $C^{1,\alpha}$ boundary.

\subsection{Proof of Theorem 1.1}

Since $\Omega$ is a bounded domain in $\CC$ with $C^{1,\alpha}$ boundary for some $\alpha>0$, there is $\rho\in C^{1,\alpha}(\CC)$
with $\Omega=\{z\in \CC: \rho(z)<0\}$ and there is a scalar constant $C_0$ such that
 $$
 {1\over C_0}\le |\nabla \rho|\le C_0, \  \hbox{ if  } \rho(z)=0.
 $$ 
 The continuity  of $|\nabla \rho|$ implies that there is an $\epsilon_0>0$ such that
$$
{1\over 2 C_0}\le |\nabla \rho(z)| \le 2C_0,\quad \hbox{ if } \ -\epsilon_0\le \rho(z)\le 0.
$$
  For $\ell\in \NN$ such that ${1\over \ell}<\epsilon_0$, define
$$
\Omega_\ell=\{z\in \CC: \rho(z)<-\ell ^{-1}\}.\leqno(5.9)
$$
Then $\d \Omega_{\ell}$ is uniformly $C^{1,\alpha} $ boundary for all $\ell\ge 2/\epsilon_0$ and
$$
\Omega_\ell\subset \Omegabar_\ell \subset \Omega_{\ell+1}\subset  \Omegabar_{\ell+1} \subset \Omega 
\quad\hbox{and}\quad
\lim_{\ell \to \infty} \Omega_{\ell} =\Omega. 
$$
Moreover, $C_{\Omega_\ell}$ is a constant depending only on $2C_0$ and $\|\rho\|_{C^{1,\alpha}(\Omegabar)}$ uniformly on $\ell> 2/\epsilon_0$.
Therefore, $C_{\Omega_\ell}=C_\Omega$ for all $\ell>2/\epsilon_0$.

Notice that
$$
f*\chi_\epsilon \in C^\infty_{(0,1)}(\Omega_\ell^n)\leqno(5.10)
$$
is $\dbar$-closed in $\Omega_\ell$ if $\epsilon<\dist(\d \Omega_\ell, \d \Omega)$/n.  By Theorem 4.3,  we have
$$
\|S_\ell[f] \|_{ L^p(\Omega_\ell^n)} \le C_\Omega \|f\|_{L^p_{(0,1)} (\Omega_\ell ^n)},\quad \hbox{for } 1\le  p\le \infty, \leqno(5.11)
$$
where $C_\Omega =C(C_0 \|\rho\|_{C^{1,\alpha}(\Omegabar)} )^n $ is a constant independent of $p$ and $\ell$. For any $1<p<\infty$, since the unit ball is weakly compact  in $L^p(\Omega_\ell)$, there is a subsequence
$\{S_{\ell_j}[f]\}_{j=1}^\infty $ converges to a function in $L^p(\Omega)$, denoted by $\tilde{S}[f]$ weakly on $L^p(\Omega_\ell)$ for any $\ell\ge 2/\epsilon_0$. Thus,
$$
\|\tilde{S}[f]\|_{L^p(\Omega_\ell^n)} \le C_\Omega \|f\|_{L^p_{(0,1)}(\Omega_\ell^n)}\le C_\Omega\|f\|_{L^p_{(0,1)}(\Omega^n)},\quad \ell\ge 2/\epsilon_0.\leqno(5.12)
$$
This implies that $\tilde{S}[f]\in L^p(\Omega^n)$ and
$$
\|\tilde{S}[f]\|_{L^p(\Omega^n)} \le C_\Omega \|f\|_{L^p_{(0,1)}(\Omega^n)}. \leqno(5.13)
$$
By the uniqueness of weak limit for each $L^p(\Omega^n)$, one has $S[f]=\tilde{S}[f]$  for all $p\in (1,\infty)$. Since $C_\Omega$ in (5.13) does not depend
 on $p$, letting $p\to\infty$, one has
$$
\|\tilde{S}[f]\|_{L^\infty(\Omega)^n} \le C_\Omega \|f\|_{L^\infty_{(0,1)} (\Omega^n)}.\leqno(5.14)
$$
Since $S_\ell [f]$ is the canonical solution for $\dbar u=f$ in $\Omega_\ell$, it is easy to check $\dbar \tilde{S}[f]=f$ in $\Omega$ in the sense of distribution.
Moreover, for any $z\in \Omega$, one has
$$
\int_{\Omega^n} \tilde{S}[f](w) K_{\Omega^n}(z,w) dv(w)=\lim_{\ell\to \infty} \int_{\Omega_\ell^n} S_\ell [f] (w) K_{\Omega^n}(z,w) dv(w)=0.\leqno(5.15)
$$
Therefore, $\tilde{S}[f]$ is the canonical solution of $\dbar u=f$ in $\Omega$. So, $S[f]=\tilde{S}[f]$. Moreover,
$\|S[f]\|_{L^\infty(\Omega^n)}\le C_\Omega \|f\|_{L^\infty_{(0,1)}(\Omega^n)}$. Therefore,   the proof of Theorem 1.1 is complete.\epf

\begin{remark} In fact, for any $\dbar$-closed $(0,1)$-form $f$ in $L^p_{(0,1)}(\Omega^n)$, we have
$$
\|S[f]\|_{L^p(\Omega)}\le C_\Omega \|f\|_{L^p_{(0,1)}(\Omega^n)},\quad \hbox{for all } 1<p\le \infty.\leqno(5.16)
$$
\end{remark}
\section{Remarks}

For any $\gamma \in [0, 1)$, we choose $\epsilon$ such that $n \epsilon=1-\gamma$. Thus, for example, $n=3$,
\begin{eqnarray*}
\lefteqn{d_{\Omega}(w)^{-\gamma}  |k(z_j, w_j)|{1\over |z_j-w_j|^2+|w_k-z_k|^2}{1\over |w_k-z_k|^2+|w_\ell-z_\ell|^2}}\\
&\le &C d_{\Omega}(w)^{-\gamma} {1\over |z_j-w_j|^{1+(n-1)\epsilon}}{|z_j-w_j|^{(n-1)\epsilon} \over |z_j-w_j|^2+|w_k-z_k|^2}{1\over |w_k-z_k|^2+|w_\ell-z_\ell|^2}\\
&\le& C{1\over d_\Omega (w)^\gamma |z_j-w_j|^{1+(n-1)\epsilon}}{1\over |z_k-w_k|^{2-\epsilon}|w_\ell-z_\ell|^{2-\epsilon}}
\end{eqnarray*}
and
$$
\int_{\Omega^3} \Big({1\over d_\Omega (w)^\gamma |z_j-w_j|^{1+(n-1)\epsilon}}{1\over |z_k-w_k|^{2-\epsilon}|w_\ell-z_\ell|^{2-\epsilon}}\Big)^{p'}dv(w)
\le C_\epsilon.
$$
for any $1<p'<{4-\epsilon \over 4-2\epsilon}$

By the arguments given in Sections 4 and  5, we have proved the following theorem.

\begin{theorem} Let $\Omega$ be a bounded domain in $\CC$  with $C^{1,\alpha}$ boundary. 
Let $f=\sum_{j=1}^n f_j d\zbar_j \in L^\infty_{(0,1)}(\Omega^n)$ be $\dbar$-closed. Then
there is a positive constant $C_\Omega$ depending only on $C^{1,\alpha}$ regularity of $\d \Omega$ such that
$$
\|S[f]\|_{L^\infty (\Omega^n) }\le {C_\Omega \over (1-\gamma)^n }\sum_{k=1}^n  \| d_\Omega (z_k)^\gamma f_k(z) \|_{L^\infty (\Omega^n)},  
$$
for any $0<\gamma<1$.
\end{theorem}

\fontsize{11}{11}\selectfont
 \bigskip
 
\noindent Department of Mathematics, University of California, Irvine, CA 92697-3875, USA

 \bigskip
\noindent Email addresses:  \ sli@math.uci.edu

\end{document}